\documentclass[12pt,a4paper]{article}
\usepackage{amsmath,amssymb,amsthm,amsfonts,amscd,euscript,verbatim}
\usepackage{hyperref}
\usepackage{newlfont}
 \theoremstyle{plain}
\newtheorem{theo}{Theorem}[section]

\newtheorem{pr}[theo]{Proposition}

 \newtheorem{lem}[theo]{Lemma}
 \newtheorem{coro}[theo]{Corollary}
  
\theoremstyle{remark}
\newtheorem{rema}[theo]{Remark}

\theoremstyle{definition}
\newtheorem{defi}[theo]{Definition}
\newtheorem*{notat}{Notation}

\newcommand\ob{^{-1}}

\newcommand\obj{Obj}
\newcommand\mo{Mor}
\newcommand\id{id}

\newcommand\z{{\mathbb{Z}}}
\newcommand\q{{\mathbb{Q}}}
\newcommand\af{\mathbb{A}}
 
\newcommand\p{\mathbb{P}}

\newcommand\pt{pt}
\newcommand\eps{\varepsilon}

\newcommand\ns{\{0\}}
\DeclareMathOperator\homm{\operatorname{Hom}}

\DeclareMathOperator\inli{\varinjlim}

\newcommand\spe{\operatorname{Spec}}

\newcommand\co{\operatorname{Cone}}

\DeclareMathOperator\cha{\operatorname{char}}

\newcommand\znz{\z/l^n\z}
\newcommand\zlz{\z/l\z}

\newcommand\com{{\mathbb{C}}}

\newcommand\modd{\operatorname{mod}}
\newcommand\tmodd{\operatorname{TopMod}}

\newcommand\znzr{\z/l^n\z(r)}
\newcommand\znrz{\z/l^n\z(r)}
\newcommand\znrpz{\z/l^n\z(r')}
\newcommand\znroz{\z/l^n\z(r-1)}
\newcommand\gi{\mathfrak{i}}
\newcommand\gj{\mathfrak{j}}
\newcommand\gff{\mathcal{F}}

\newcommand\zol{{\mathbb{Z}[\frac{1}{l}]}}
\newcommand\op{{}^{-1}}

\newcommand\ds{D^b_c\,Sh^{et}_{\znz-\modd}} 
\newcommand\dsx{D^b_c\,Sh^{et}_{\znz-\modd}(X)} 
\newcommand\dsy{D^b_c\,Sh^{et}_{\znz-\modd}(Y)} 

\begin{document}

 \title
 {A 'fat hyperplane section' weak Lefschetz (in arbitrary characteristic) and Barth-type theorems}
 \author{M.V. Bondarko
   \thanks{ The work is supported by RFBR
(grants no. 10-01-00287, 11-01-00588, and 12-01-33057), by the Saint-Petersburg State University research grant no. 6.38.75.2011, and by the Federal Targeted Programme "Scientific and Scientific-Pedagogical Personnel of the Innovative Russia in 2009-2013" (Contract No. 2010-1.1-111-128-033). }}
 \maketitle
\begin{abstract}

In this  paper we prove a certain 'algebraic fat hyperplane section'  Weak Lefschetz-type result for $\znz$-\'etale cohomology (of locally a complete intersection variety $X'$ that is not necessarily proper); somewhat similar statements for 
$X'\stackrel{s}{\to} \com\p^N$ were previously established by Goresky and MacPherson using stratified Morse theory. In contrast to their 'topological' statements (that relate  $\pi_i(X')$ with $\pi_i(s\ob(\com\p^{N-b}_{\eps}))$, where $\com\p^{N-b}_{\eps}$ is the $\eps$-neighbourhood of $\com\p^{N-b}$ in $\com\p^N$ for a small $\eps>0$) 
 we  formulate and easily prove our (cohomological) result using purely sheaf-theoretic methods; this makes it independent from the base field characteristic. Our proof is quite short; we apply an argument similar 
 to the one used by Beilinson in order to establish a Weak Lefschetz theorem for general hyperplane sections of a smooth $X'\subset \p^N$. 
 Considering our main theorem as a substitute of the one of Goresky and MacPherson seems to be an important idea; it allows to combine the implications of the fact mentioned (that were previously known for complex varieties  and quasi-finite $s$ only) 
 with  (very convenient) sheaf-theoretic methods (in particular, with proper and smooth base change).
So we extend (generalizations of) 
the seminal results of Barth and others to base fields of arbitrary  characteristics ($\neq l$; we even obtain  certain similar statements over an arbitrary $S/\spe \zol$). We obtain certain statements 
concerning the lower cohomology 
for  closed subvarieties of $\p^N$ of small codimensions and for their preimages with respect to proper morphisms (that are not necessarily finite); to this end we apply the methods developed by Deligne, Fulton, and Lazarsfeld.
Note: whereas a certain Barth-type theorem for subvarieties of $\p^N$ (over any field $K$) 
was proved by Lyubeznik, he proved nothing about their preimages.



\end{abstract}

\tableofcontents

 \section*{Introduction}

In contemporary algebraic geometry there are a lot of 'descendants' of the classical Weak Lefschetz theorem (a rich collection of those, that includes the theorem of Barth on the cohomology of closed subvarieties of $\p^N$ of low codimensions, can be found in \cite{lazar}). Whereas for the 'ordinary' Weak Lefschetz theorem the Artin's vanishing yields a 'purely algebraic' proof (that works over an arbitrary characteristic base field $K$),  the proofs of several other results are often somewhat 'topological', or use Hodge theory, or De Rham cohomology. These methods of the proofs restrict the facts obtained to the 
 case $\cha K=0$ (the corresponding literature is so vast that the author will not attempt to mention all of these results).

In particular, one of the powerful tools for studying these Weak Lefschetz-type questions for $K=\com$  (as shown in particular in \S9 of \cite{fulaz}; cf. also \S3.5.B of \cite{lazar}; the corresponding method was proposed by Deligne) is the 'fat (multiple) hyperplane section' weak Lefschetz theorem formulated in \S II1.2 of \cite{gorma}. For a quasi-finite morphism $s:X'\to \com\p^N$ ($N>0$, $X'$ is  locally a complete intersection variety) and any  small enough  $\varepsilon>0$  it states that the lower homotopy groups of $X'$ are isomorphic to those of $s\ob (\com\p^{N-b}_\varepsilon)$, where  $\com\p^{N-b}_\varepsilon$ is the $\varepsilon$-neighbourhood (in the sense of some Riemannian metric for $\com\p^N$) for the standard embedding $\gi$ of $\com\p^{N-b}$ into $\com\p^N$ ($0<b<N$) (cf. the caution below).  Note that even the formulation of this statement requires some 'topology', whereas the proof given in ibid. heavily relies on  Stratified Morse Theory. Thus, both the formulation and the proof of loc. cit. are far from being 'algebraic' (and in the opinion of the author, the proof mentioned is  extremely complicated, and is very long to write down in full detail).

The current paper grew our from the following simple observation: the direct limit of the  $\znz$-cohomology of  $s\ob (\com\p^{N-b}_\varepsilon)$ when $\varepsilon\to 0$ 
is just the (hyper)cohomology of $\com\p^{N-b}$ with coefficients in $R\gi^*Rs_{*}\znz_{X'}$ (see Remark \ref{rara}(\ref{i3}) below). 
Then it was a simple exercise to verify that the lower $\znz$-cohomology of $X'$ is isomorphic to the one of $R\gi^*Rs_{*}\znz_{X'}$; this is true both for singular and for \'etale cohomology. 
Though our proof relies on certain results of \cite{bbd} and \cite{kw} (since the perverse $t$-structure is a very convenient tool for our purposes), its base is just Artin's vanishing. 
 Certainly, the fact obtained (in the setting of \'etale cohomology) does not depend on the choice of  
 $K$ (and $p$); this is one of its major advantages over the results of \cite{gorma}. To the opinion of the author our proof  below is much easier than the corresponding one of Goresky and MacPherson (though it is somewhat technical). Though our reasoning has nothing to do 
  with the arguments of ibid., 
  we should note that somewhat similar methods were previously used by Beilinson in order to establish a Weak Lefschetz-type statement for general hyperplane sections of a smooth $X'\subset \p^N$ (see Lemma 3.3 of \cite{bei87} and Remark \ref{rhens}(3) below). Yet considering our main theorem as a substitute of the 
results of \cite{gorma} seems to be an important idea that  has several nice consequences (some of them are really new, especially in the case $p>0$). 
In particular, we immediately obtain a (generalized) Weak Lefschetz theorem for projective ('almost') locally set-theoretic complete intersections in arbitrary characteristic ($\neq l$). 

A disadvantage of our methods is that 
for $K=\com$ they yield no information on   (the cohomology of) $s\ob (\com\p^{N-b}_\varepsilon)$ for any particular $\varepsilon>0$. 
On the other hand, we do not demand $s$ to be quasi-finite (see also Remark \ref{rara} 
 below for a further discussion on the comparison of our results with those of \cite{gorma}, including Theorem II1.1 of ibid.). Besides, our 'sheaf-theoretic' statement can be easily combined with  proper and smooth base change; see (the proofs of) Proposition \ref{pinsub} and Corollary \ref{cpro} below. This allows us to use it as a substitute of the theorem of \cite{gorma} in the argument of Deligne described in \S9 of \cite{fulaz} (see also  \S3.5.B of \cite{lazar}). As a consequence, we
obtain certain cohomological analogues of  the  results of loc. cit.  over arbitrary ($\neq l$) characteristic fields.  
In particular, we easily extend 
the theorems of Barth %
on the lower cohomology  
for  
closed subvarieties of $\p^N$ of small codimensions and for their preimages with respect to proper morphisms (we do not require those to be finite in contrast to Corollary 9.8 of \cite{fulaz}) to the case of arbitrary characteristic ($\neq l$).
More generally, we compute the lower cohomology of the preimages of the diagonal in $(\p^N)^q$ with respect to proper morphisms. 
We even deduce  certain versions of these statements over an arbitrary $S/\spe \zol$ (they are formulated in terms of $S$-sheaf \'etale cohomology of $S$-schemes).

Taking all of this into account, the author hopes that his methods will become a useful tool for studying (various) Weak Lefschetz-type questions (at least, 'higher degree' ones). We should certainly mention  here that Lyubeznik has also applied certain 'sheaf-theoretic' methods to the study of the lower cohomology of subvarieties of $\p^N$ for an arbitrary $p$; in Theorem 10.5 of \cite{lyubez} he proved several results that are (somewhat) stronger than our Theorem \ref{tbarth}(II1) (cf. also Remark \ref{rara}(\ref{i4}) below). 
 Yet it seems that the methods of ibid.
(that are very interesting and quite distinct from our ones) cannot say much on the cohomology of the preimages of subvarieties in $\p^N$ (and in $(\p^N)^q$, with respect to  proper morphisms; cf. Remark \ref{rlyu}(5)).

{\bf A caution}: for $K=\com$ the natural morphism  $h^i:H^i(X',\znz)\to H^i(s\ob (\com  \p^{N-b}),\znz)$ certainly factorizes as the composition $H^i(X')\stackrel{f^i}{\to} \inli_{\eps\to 0}H^i(s\ob (\com\p^{N-b}_\varepsilon))\stackrel{g^i}{\to}H^i(s\ob (\com  \p^{N-b}),\znz)$
(for any $ i\ge 0$). In the proof of Corollary \ref{cpro}(1) below we describe an algebraic analogue of this factorization (through $H^i(\p^{N-b},R\gi^*Rs_{*}\znz_{X'})$; unfortunately, the proof is somewhat formal). It turns out that this factorization simplifies the study of weak Lefschetz-type questions. In particular,  the results of \S II1.2 of \cite{gorma} (in the case $K=\com$) and  our Theorem \ref{tgara} (for the general case) yield that $f^i$ is bijective for 'small enough' $i$  (even if $s$ is not proper). 
Note that this result does not extend to $h^i$ (in the general case); in particular,  $s\ob (\com \p^{N-b})$ can be empty. Yet if $s$ is proper over 
some open $U\subset X$, $Z\subset O$,  then $g^i$ is necessarily an isomorphism for all $i$ (our proof of this result is a more or less easy combination of smooth and proper base change theorems). 
Thus one may think about $H^*(\p^{N-b},R\gi^*Rs_{*}\znz_{X'})$ as of an 'approximation' to the cohomology of $s\ob (\p^{N-b})$ that has several nice properties. Our Theorem \ref{tgara} yields that $f^i$ is an isomorphism for lower $i$ (one may say that this is 'a Goresky-MacPherson-type' result); in Corollary \ref{cpro}(1) we use it in order to establish a ('true') Weak-Lefschetz-type statement, which we actively apply in \S\ref{sappl}.

Lastly we note that our methods also yield (without any problem at all) a certain generalization of our basic result to the relative setting (for example, to schemes over $\spe\zol$ instead of $K$-varieties). Yet this expansion would be stated in terms of perverse sheaves (over the base), and at the moment the author knows no 'visualizable' applications for it (still see Remark \ref{rara}(\ref{i6}) below).

Now we describe the contents of the paper. First we give a 'brief plan' of it. The main result of the paper is Theorem \ref{tgara} (our 'fat hyperplane section Weak Lefschetz'); we use it along with a simple (smooth and proper) base change argument in order to deduce Corollary \ref{cpro} (that states that under certain restrictions a 'true Weak Lefschetz-type statement' holds). We apply 
  part 1 of the Corollary to certain $G_m^s$-bundles (constructed by Deligne) in Theorem \ref{tbarth}(I); 
then we deduce a certain
Barth-type theorem (in Theorem \ref{tbarth}(II); see Remark \ref{rlyu}(6) for a generalization of this result to the case of schemes over an arbitrary  $S/\spe \zol$). 

In \S\ref{sprel} we recall some basics on the derived categories of (constructible) $\znz$-sheaves, on functors between them, and on the perverse $t$-structure. 

In \S\ref{sbase} we prove our 'fat multiple hyperplane section' weak Lefschetz-type Theorem \ref{tgara}. We also make several remarks
(on possible modifications of loc. cit.; a reader that is only interested in the results of \S\ref{sappl} may skip the rather long Remark \ref{rara}), and show how one can use proper and smooth base change in order to calculate the (hyper)cohomology of $R\gi^*Rs_{*}\znzr_{X'}$.

In \S\ref{sappl} we describe some applications of Theorem \ref{tgara}. 
Following \S9 of \cite{fulaz} (see also \S3.5.B of \cite{lazar}), we consider certain (locally trivial) $G_m^s$-bundles. This  relates the cohomology of a (certain $G_m^{q-1}$-bundle over) a variety $Y$ that is proper over $(\p^N)^q$ with the one of  the preimage of the diagonal.
As a result, we extend the corresponding theorems  of Barth  to the case of an arbitrary $p\neq l$ (and to $\znz$-coefficients; we also obtain the relative versions mentioned above).

In \S\ref{srem} we make several  remarks on the calculating of the (hyper)cohomology of 
$R\gi^*(-)$  via henselizations (that are related with author's ideas on the proof of a certain Weak Lefschetz-type statement for torsion motivic cohomology). 

The author is deeply grateful  prof. H. Esnault, to prof. M. Hoyois,   to prof. G. Lyubeznik, to prof. I. Panin, and to  prof. A. Vistoli 
 for their  interesting remarks. He would also like to express his gratitude
 prof. M. Levine and to the Essen University, as well as to prof. F. D\'eglise and to Unit\'e de math\'ematiques pures et appliqu\'ees of
\'Ecole normale sup\'erieure de Lyon for the wonderful working conditions (in the summer of 2012 and January of 2015, respectfully).

 \begin{notat}


$K$ will be our  base field of characteristic $p$ ($p$ can be zero). 
The reader may (always) assume 
that  $K$ is algebraically closed; in particular, then one can ignore some of the 
statements from \S\ref{sprel}. 
 We will denote by $X_{red}$ the reduced scheme associated to a scheme $X/K$. 

 $\pt$ is a point, $\af^N$ is the $N$-dimensional
affine space (over $K$), 
$\p^N$ is the projective space of dimension $N$, $G_m=\af^1\setminus \ns$ is the multiplicative group scheme.

We will say that a variety has dimension $m$ only in the case when it is equidimensional of dimension $m$. 
We will say that a morphism of varieties has relative dimension $\le e$ if the dimension of all of its fibres is $\le e$.

We will call a equidimensional  variety $V$ locally a set-theoretic complete intersection (or just an LSTCI) if it possesses a Zariski cover each of whose  components is isomorphic to an open subvariety of a  set-theoretic complete intersection in $\p^N$ (for some $N\ge 0$). This is certainly equivalent to the existence of a Zariski cover each of whose  components is
isomorphic to  a subvariety of $\p^N$ (or of $\af^N$) 
whose closure is defined by $N-d$ equations, where $d$ is the dimension of $V$ (cf. \cite{lyubez}). 

Throughout the paper we will be interested in \'etale cohomology with $\z/l^n\z(r)$-coefficients, where $l$ is a prime $\neq p$, $n>0$, and $r\in \z$. 
One can assume $l,n$ to be fixed; 
besides, the case $r=0$ is already interesting enough (and is equivalent to the general case if $K$ is an algebraically closed). Moreover,   one may also consider cohomology with coefficients in arbitrary $l$-torsion locally constant \'etale sheaves (cf. Remark \ref{rara}(\ref{i1})).

Below we will consider a morphism $s_X:X'\to X$ and a closed embedding $\gi:Z\to X$; in the introduction and in the abstract we took $X=\p^N$, $Z=\p^{N-b}$, and wrote $s$ instead of $s_X$. Besides, following \cite{bbd} we will always omit $R$ in $Rf_*,\  Rf_!,\ Rf^*$, and $Rf^!$  (for $f$ being a finite type morphism of varieties); also, we will ignore the difference between cohomology (of sheaves) and hypercohomology (of complexes of sheaves).  $t$ will denote the (self-dual) perverse $t$-structure for $\ds(-)$; see below.


$\ds(Y)$ will denote the derived category of complexes of \'etale $\znz$-module sheaves over a variety $Y$ with 
constructible cohomology. We will introduce some (more) notation for these categories in Proposition \ref{pds}. 

For an (additive) category $C$ we will denote by $C(A,B)$ the group $\mo_C(A,B)$.

\end{notat}

\section{Some preliminaries on derived categories of sheaves and the perverse $t$-structure for them}\label{sprel}

In this section we introduce some notation and remind some of the properties of
the categories  $\ds(-)$ and functors between them; we relate 
 those with the  properties of the 'usual' \'etale (hyper)cohomology of complexes of sheaves. 
We also describe the properties of the perverse $t$-structure for $\ds(-)$.
The statements of this section are well-known (possibly except Proposition \ref{pbch}(I1), which is  quite easy). Most of them were (essentially) proved in  SGA4, SGA41/2, or in \cite{bbd}. 

We do not mention the \'etale fundamental group here; though we 
consider it in Corollary \ref{cpro}(2,3) below, 
 these statements do not seem to be 
 important (to the author).


\begin{pr}\label{pds}

Let $f:X\to Y$ be a morphism of $K$-varieties.

\begin{enumerate}


\item\label{icoh} For $C\in \obj D(Sh^{et}(X,\znz-\modd)$, $i\in\z$,  we will denote
$D(Sh^{et}(X,\znz-\modd))(\znz_X,C[i])$ 
  by $H^i(X,C)$; 
this is the 'usual' $i$-th (hyper)cohomology of $X$ with coefficients in $C$.


\item\label{idscat} For any variety $X/K$ there is a full 
 triangulated subcategory $\dsx\subset D(Sh^{et}(X,\znz-\modd))$. 

For $C\in \obj D(Sh^{et}(X,\znz-\modd)),\ s\in \z$ we will denote $C\otimes \znz(s)$ by $C(s)$ (whereas $C[i]$ for $i\in \z$  denotes the shift of $C$ by $i$ 'to the left'). We have that $\dsx(s)=\dsx$ for any $s\in \z$.

\item\label{ifun}  
The following pairs of adjoint functors
 are defined: $f^*: D(Sh^{et}(Y,\znz-\modd) \leftrightarrows D(Sh^{et}(X,\znz-\modd):f_*$ and $f_!: D(Sh^{et}(X,\znz-\modd) \leftrightarrows D(Sh^{et}(Y,\znz-\modd):f^!$; they restrict to functors 
 $f^*: \ds(Y) \leftrightarrows \dsx:f_*$ and $f_!: \dsx \leftrightarrows \dsy:f^!$.
Any of these four types of functors (when $f$ varies) yields a  2-functor from the category of 
$K$-varieties to the 2-category of triangulated categories. 

\item\label{ik} $D(Sh^{et}(\spe K,\znz-\modd))$ is isomorphic to the derived category $ D(\znz[G]-\tmodd)$ of topological $G$-modules, where $G$ is the absolute Galois group of $K$ (endowed with the corresponding topology); 
this is an isomorphism of tensor triangulated categories. 

\item\label{icohh}
Consider 
 the structure morphisms $x:X\to \spe K$ and $y:Y\to \spe K$. Then for any $C\in \obj \ds(Y)$ we have that
 $H^*(Y,C)\cong H^*(\spe K, y_*C)$. 
 Moreover,  the 'usual' morphisms of cohomology $H^*(Y,C)\to H^*(X,f^*C)\cong H^*(\spe K, x_*f^*C)=H^*(\spe K, y_*f_*f^*C)\cong H^*(Y, f_*f^*C)$ come from 
 the adjunction $f^*: \ds(Y) \leftrightarrows \ds(X):f_*$.

\item\label{iupstar}  $f^*$ is symmetric monoidal; $f^*(\znz(s)_Y)=\znz(s)_X$ for any $s\in \z$. 

\item \label{ipur}
$f_*= f_!$ if $f$ is proper; 
if $f$ is an open immersion, 
then $f^!=f^*$.

\item\label{ikun} Let $K$ be 
algebraically closed. Then
for any $s\in \z$, $X/K$, one can choose an isomorphism $\znz_X\cong \znz(s)_X$ that is functorial in $X$. Besides,
for 
 $Z=X\times Y$  we have the Kunneth formula for the cohomology of $Z$, i.e., 
for the corresponding structure morphisms $x,y,z$ we have that
$z_*\znz_Z\cong x_*\znz_X \otimes y_*\znz_Y$. 

\item\label{iglu}
If $i:Z\to X$ is a closed immersion, $U=X\setminus Z$, $j:U\to X$ is the complementary open immersion, then
the pairs of morphisms (that come from the adjunctions of assertion \ref{ifun})
\begin{equation}\label{eglu}
j_!j^*(M) \to M
\to i_*i^*(M)
\end{equation}
  and 
\begin{equation}\label{eglu1} i_*i^!M \to M \to j_*j^*M \end{equation} can be completed to
distinguished triangles.

\item\label{igysin} 
In the setting of the previous assertion, if $i$ is a closed immersion of smooth varieties everywhere of some codimension $s>0$, $U=X\setminus Z$, $j:U\to X$ is the corresponding embedding, then for any $C\in \obj D(Sh^{et}(X,\znz-\modd))$ there exists a distinguished triangle \begin{equation}\label{eglu2} i_*i^*C(-s)[-2s]\to C  \to j_*j^*C.\end{equation} 

\end{enumerate}
\end{pr}
\begin{proof} The statements are well-known (see  SGA4 and SGA41/2 for the proofs) and most of them were (actively) used in \cite{bbd} (see also Theorem 6.3 of \cite{eke}). 
We will only give a few 
(more precise) references.

Assertion \ref{idscat} is mostly given by Corollary 1.5 of \cite{sga45fin}.

Assertion \ref{icohh} can be deduced from the 'classical' properties of cohomology by applying adjunctions.

Assertion \ref{ikun} 
follows easily from Corollary 1.11 of loc. cit.


Assertion \ref{igysin} can be obtained from the distinguished triangle (\ref{eglu1}) by applying purity; cf.  \S3.2 of \cite{sga418}.

\end{proof}

We will also recall certain properties of  base change transformations (closely following \S4 of \cite{sga417}).

\begin{defi}\label{dbch}
For a commutative square 
\begin{equation}\label{ebch} \begin{CD} 
 X' @<{\gi'}<< Z'\\
@VV{s_X}V @VV{s_Z}V\\
X @<{\gi}<< Z
\end{CD}\end{equation}
of  morphisms of varieties 
 and $C\in \obj\ds(X')$ we call the composition of the morphisms $\gi^*s_{X*}C\to   \gi^* s_{X*}\gi'_*\gi'^*C=  \gi^*\gi_* s_{Z*} \gi'^*C       \to s_{Z*}\gi'^*C$ coming from the corresponding adjunctions the {\it base change} morphism; we denote the corresponding  natural transformation $\gi^*s_{X*}\implies s_{Z*}\gi'^*$  by $B(s_X,\gi)$.
\end{defi}

  \begin{pr}\label{pbch}
  
I In the setting of   Definition \ref{dbch} the following statements are fulfilled.
  
  1. The morphism $s_{X*}C\to  s_{X*}\gi'_*\gi'^*C$ that comes from the adjunction $\gi'^*: \ds(X') \leftrightarrows \ds(Z'):\gi'_*$  equals the composition of $\gi_*(B(s_X,\gi)(C))$ with the morphism $s_{X*}C\to \gi_*\gi^*s_{X*}C$ coming from the adjunction  $\gi^*: \ds(X) \leftrightarrows \ds(Z):\gi_*$.
  
  2.  Let (\ref{ebch}) be a cartesian square. Then $B(s_X,\gi)$ is an isomorphism either  if $\gi$ is smooth 
   (this statement is called {\it smooth base change}), or if $s_X$ is proper (this is {\it proper base change}).

II 
Base change transformations respect compositions, i.e., for a commutative diagram
$$ \begin{CD} Z'@>{i_1'}>> O' @ >{i_2'}>> X'\\
@VV{s_Z}V
@VV{s_O}V @VV{s_X}V\\
Z@>{i_1}>> O @ >{i_2}>> X \end{CD}$$
we have that $B(s_X,i_2\circ i_1)=B(s_O,i_1 )(i'_{2*}(-)) \circ i_1^*(B(s_X,i_2 ) ) $.

\end{pr}
\begin{proof}
I1. So, we have the composition $s_{X*}C\to \gi_* \gi^*s_{X*}C\to   \gi_*\gi^* s_{X*}\gi'_*\gi'^*C= \gi_* \gi^*\gi_* s_{Z*} \gi'^*C       \to \gi_* s_{Z*}\gi'^*C=s_{X*}\gi'_* \gi'^*C$; hence it suffices to verify that the composition  $s_{X*}\gi'_* \gi'^*C=\gi_* s_{Z*} \gi'^*C\to
 \gi_* \gi^*\gi_* s_{Z*} \gi'^*C       \to \gi_* s_{Z*}\gi'^*C$ is the identity.
The latter statement is immediate from the fact the the composition of transformations $\gi_*\implies \gi_* \gi^*\gi_*\implies \gi_*$
(coming from the  adjunction $\gi^*: \ds(X) \leftrightarrows \ds(Z):\gi_*$) is the identical transformation of $\gi_*$ (this is a basic general property of adjunctions).

2. These are the derived versions of the classical base change isomorphism statements; 
see Corollary 1.2 of \cite{sga416} and Theorem 4.3.1 of \cite{sga417}, respectively. 

II This is also a well-known statement (see \S2 of ibid.).
We should check that $B(s_X,i_2\circ i_1)$ equals
the composition $(i_2\circ i_1)^*s_{X*}\implies    (i_2\circ i_1)^*   s_{X*}i'_{2*}i'^*_2=     i^*_{1}    i^*_{2}  i_{2*}s_{O*} i'^*_2\implies
             i^*_{1}   s_{O*} i'^*_{2} \implies   i^*_{1}   s_{O*}  i'_{1*}i'^*_1 i'^*_{2} = i^*_{1}   i_{1*}   i'^*_{2}s_{Z*}  i'^*_{1}    \implies   s_{Z*}  i'^*_{1}  i'^*_{2} $; this is obvious. 

\end{proof}


Now we recall some of the properties of the  perverse $t$-structure $t$ for $\ds(-)$ (that corresponds to the self-dual perversity denoted by $p_{1/2}$ in \cite{bbd}).
We will not need much of them (we even will not need a definition for $t$); so we just recall some of the results of \cite{bbd} and verify briefly that two statements from \cite{kw} can be carried over to our setting (of $\znz$-module sheaves for an arbitrary $K$). 

\begin{pr}\label{ppt}
For the perverse $t$-structure $t$ (given by the couple ($\ds(-)^{p \le 0},\break \ds(-)^{p \ge 0}$); cf. ibid.) corresponding to the perversity $p_{1/2}$ the following statements are fulfilled.
\begin{enumerate}

\item\label{ifield} If $X=\spe K$, $K$ is a field, then $t$ is just the canonical $t$-structure for $\ds(\spe K)$ (that is compatible with the canonical $t$-structure for $D(\znz[G]-\tmodd)$, where $G$ is the absolute Galois group of $K$).
In particular, if $C\in \obj\ds(\spe K)^{p \ge 0}$ then $H^i(\spe K,C)=0$ for any $i<0$.

\item\label{ihsss} Let $L$ be the algebraic closure of a field $K$, $f:X\to Y$  is a morphism of $K$-varieties; denote by 
$f_{L,red}:X_{L,red}\to Y_{L,red}$ the corresponding morphism of $L$-varieties.

Suppose that for some $N,r\in \z$ we have that $H^i(-,\znrz)(f_{L,red})$ is bijective for all $i<N$, and is injective for $i=N$.
Then $H^i(-,\znrz)(f)$ is also bijective for all $i<N$, and is injective for $i=N$.


\item\label{itaff} Suppose that a variety $U$ can be presented as the union of $b$ open affine subvarieties (for some $b$); $u:U\to \spe K$ is its structure morphism. Then $u_!(\ds(U)^{p\ge 0})\subset \ds(K)^{p\ge 1-b}$.

\item\label{itexdimd} Let $f:X\to Y$ be a morphism of relative dimension $\le d$. Then $f_*(\ds(X)^{p\ge 0})\subset \ds(Y)^{p\ge -d}$.

\item\label{ilstci} Let $X$ be  locally a set-theoretic complete intersection (in the sense of the Notation, i.e., suppose that it possesses a Zariski cover each of whose  components is isomorphic to an open subvariety of a  set-theoretic complete intersection in $\p^N$ for some $N\ge 0$) of dimension $a$. Then $\znz(r)_X\in \ds(X)^{p\ge a}$ for any $r\in \z$. 

\end{enumerate}
\end{pr}
\begin{proof}

Most of these statements were 
proved in \cite{bbd} (note that the proofs work in the case of $\znz$-coefficients; see \S4.0 of \cite{bbd}.

\begin{enumerate}

\item Immediate from the definition of $t$ (see \S2.2.12-13 of ibid.). 

\item Applying Corollary of \cite{sga48}, we can pass to the case of a perfect $K$  (so, $f_{L,red}:X_{L,red}\to Y_{L,red}$ equals $f_L:X_L\to Y_L$).

We apply Proposition \ref{pds}(\ref{icohh}); denote the corresponding structure morphisms $X,Y\to \spe K$ and $X_L,Y_L\to \spe L$ by $x,y,x_L,y_L$, respectively. 
 We obtain that $H^i(\spe L, C_L)=0$ for $i< N$, where $C_L=\co(y_{L*}y_L^*\znrz_L\to y_{L*}f_{L*}f_L^*y_L^*\znrz_L)$. Obviously, this is equivalent to $C_L\in \ds(\spe L)^{p\ge N}$. Consider $C_K=\co(y_{*}y^*\znrz_K\to y_{*}f_{*}f^*y^*\znrz_K)\in D(\znz[G]-\tmodd)$.

We note that  $C_L$ can be obtained from $C_K$ via the forgetful functor $D(\znz[G]-\tmodd)\to D(\znz-\modd)$ (this is a consequence of smooth base change; note 
that one can 'pass to the limit' in Proposition \ref{pbch}(I2); see Theorem 1.1 of \cite{sga416}). Hence $C_K\in  D(\znz[G]-\tmodd)^{p\ge i}$; so it belongs to $\ds(\spe K)^{p\ge N}$. It remains to apply the previous assertion.

\item Immediate from  \S4.2.3 of \cite{bbd}.

\item The statement is contained in \S4.2.4 of ibid.

 \item 
 This is 'the easier half of' the $\znz$-coefficient analogue of Lemma III6.5 of \cite{kw}; 
 so we will only sketch the reduction of the statement to the results of \cite{bbd}. First we note that  the property of belonging to  $\ds(-)^{p\ge a}$ is Zariski-local (it is even \'etale-local immediately from the definition of $p$); 
hence we can assume that $X$ is a set-theoretic complete intersection in an open affine subvariety $V$ of 
$\p^N$ for some $N>0$. We denote by $i$ the corresponding embedding.

Since $i_*$ is conservative and $t$-exact (see \S4.2.4 of \cite{bbd}), it suffices to verify that $i_*\znz_{X} \cong i_*i^*\znz_{V} \in \ds(\p^N)^{p\ge a}$.
Now we note that $\znz_{V}\in \ds(V)^{p=N}$ (since 
$v$ is smooth, and so for $v:V\to \spe K$ the functor $v^*[N]$ 
 is $t$-exact; see 
 Proposition 4.2.5 of ibid.). Next,   for $j$ being the embedding of $U=V\setminus X\to V$ we have that $j_!(\ds(U)^{p\ge 0})\subset 
 \ds(V)^{p\ge 1+a-N}$; see   \S4.2.3 of ibid (note that $X\cap V$ is a union of $N-a$ smooth affines). Consider the following (rotation of the) distinguished triangle given 
by (\ref{eglu}) : 
$ \znz_{V} \to  i_*i^*\znz_{U} \to j_!j^*\znz_{V}[1]$. We have that $j_!j^*\znz_{V}[1]\cong j_!\znz_{V}[1]\in \linebreak
\ds({V})^{p \ge a}$; since $\znz_{V}\in \ds({V})^{p \ge a}$, the same is true for $i_*\znz_{X}$ and we obtain the result in question.

hence we can assume that $X$ is a set-theoretic complete intersection in $\p^N$ for some $N>0$; we denote by $i$ the corresponding embedding.


\end{enumerate}
\end{proof}

\begin{rema}\label{ralt} Using devisage, we could have reduced some of the results of this section to the setting of sheaves of $\zlz$-vector spaces. 
\end{rema}

\section{Our 'fat hyperplane section' weak Lefschetz-type theorem}\label{sbase}


\begin{lem}\label{lara}
Let 
 $s_X:X'\to X$ be
a   morphism of varieties; let $\gi :Z\to X$ be a closed immersion.
Denote by $x:X\to \spe K$ the structure morphism of $X$.

Then we have the following: $H^*(X',\znrz_{X'}) \cong H^*(X, s_{X*}\znzr_{X'})$, and $H^*(Z,\gi^*s_{X*}\znrz_{X'})\cong H^*(X, \gi_*\gi^* s_{X*}\znzr_{X'})$.

\end{lem}
\begin{proof}
Denote by $x'$ and $z$ the structure morphisms of $X'$ and of $Z$, respectively.

We apply Proposition \ref{pds}(\ref{icohh}). We have that $H^*(X', \znzr_{X'})\cong H^*(\spe K, x'_*\znzr_{X'})=H^*(\spe K, x_*s_{X*}\znzr_{X'})\cong H^*(X, s_{X*}\znzr_{X'})$, and $H^*(Z,\gi^*s_{X*}\znrz_{X'})\cong H^*(X, \gi_*\gi^*s_{X*}\znrz_{X'})$.

\end{proof}

Consider the morphism $\gff_r: s_{X*}\znzr_{X'}\to \gi_*\gi^* s_{X*}\znzr_{X'}$ coming from the adjunction $\gi^*: \ds(X) \leftrightarrows \ds(Z):\gi_*$. 
The base of the results of this paper is the following statement.

\begin{theo}
\label{tgara}
 
In the setting of the previous lemma, suppose that  
 $s_X$ is of relative dimension $\le e$ (i.e., the dimension of all of its fibres is $\le e$), 
$X$ is proper; for $U=X\setminus Z$
assume that  $U'=s_X\ob(U)$ is 
 locally a set-theoretic complete intersection (everywhere) of dimension $a$ (see the Notation section), and that $U$ can be presented as the union of $b$ open affine subvarieties (for some $a,b,e\ge 0$).

Then $H^i(X,-)(\gff_r): H^i(X',\znrz_{X'})\to H^i(Z,\gi^*s_{X*}\znrz_{X'})$ (see the previous lemma) 
is a bijection for $i<a-e-b$, and is an injection 
for $i=a-e-b$. 

\end{theo}
\begin{proof}

 We have the following commutative diagram  
$$\begin{CD} 
U'@>{\gj'}>> X' @<{\gi'}<< Z'\\
@VV{s_U}V
@VV{s_X}V @VV{s_Z}V\\
U@>{\gj}>>X @<{\gi}<< Z
\end{CD}$$
(both of the corner squares are cartesian); we denote the corresponding structure morphisms (whose target is $\spe K$) by $u,x,z,u',x'$, and $z'$, respectively. 

Applying (\ref{eglu}), we obtain that there exists a distinguished triangle $\gj_{!} \gj^*s_{X*}\znzr_{X'}  \to
s_{X*} \znzr_{X'}\to \gi_{*} \gi^*s_{X*}\znzr_{X'}$. 

 Proposition \ref{pds}(\ref{icohh}) 
  yields:
 we should prove that 
 the cohomology of $\spe K$ with coefficients in $x_*\gj_{!} \gj^*s_{X*}\znzr_{X'} $ vanishes in degrees  $\le a-e-b$.
To this end it suffices to verify that  $x_*\gj_{!} \gj^*s_{X*}\znzr_{X'}\in \ds(\spe K)^{p\ge a-e-b+1}$ (see Proposition \ref{ppt}(\ref{ifield})).

By smooth base change (see Proposition \ref{pbch}(I2)), 
the latter term is isomorphic to $x_*\gj_{!} s_{U*}\gj'^* \znzr_{X'}$. We transform this further into
$x_*\gj_{!} s_{U*} \znzr_{U'}=u_{!} s_{U*} \znzr_{U'}$
 (we use the properness of $x$; see Proposition \ref{pds}(\ref{ipur})). 

 We have that $\znzr_{U'}\in \ds(U')^{p\ge a}$ 
 (see Proposition \ref{ppt}(\ref{ilstci})). 
  Next,  
 $s_{U*}$  sends $\ds(U')^{p \ge a}$ into $\ds(U)^{p\ge a-e}$  (since $s_U$ is   
 of  relative dimension $\le e$; see part \ref{itexdimd} of loc. cit.). 

Now, $u_{!}$ sends  $\ds(U)^{p\ge a-e}$  into $\ds(\spe K)^{p\ge a-e-b+1}$ by part \ref{itaff} of loc. cit.
Hence  part \ref{ifield} of loc. cit. yields the result.

\end{proof}

\begin{rema}\label{rara}
 A reader that is only interested in the applications of our result described in \S\ref{sappl} may skip this remark.
\begin{enumerate}
    \item \label{i1} 
    Instead of a general $\gi:Z\to X$ one can consider $\gi$ being the natural embedding of $\p^{N-b}$ into $\p^N$; this doesn't decrease the generality of the result significantly.

Besides, instead of $\znrz$ we could have considered any locally constant $l$-torsion sheaf on $X'$ (cf. Theorem 10.5 of \cite{lyubez}). Indeed (as shown by our proof) it suffices to verify  for any locally constant $\znz$-sheaf $C/X'$ that  $\gj'^*C\in \ds(U')^{p=a}$. Now,  the latter statement can be verified \'etale-locally with respect to $U'$; so we can assume that  $\gj'^*C$ is constant and extend
 (the proof of) Proposition \ref{ppt}(\ref{ilstci}) to this setting. 


\item \label{i2} The theorem was inspired by the 
'fat multiple hyperplane section' version of Weak Lefschetz proved by M. Goresky and R. MacPherson (for complex varieties; see the theorem in \S II.1.2  of \cite{gorma}). 
 For $\gi$ being the embedding of $Z=\com \p^{N-b}$ into $X=\com\p^N$, a quasi-finite $s_X$ (i.e., $e=0$), and any small enough $\varepsilon>0$ loc. cit. states the following: for  
  the $\varepsilon$-neighbourhood (in the sense of some Riemannian metric for $X$) $Z_\varepsilon$ of $Z$ in $X$ the natural map $\pi_i(s_X\ob(Z_\varepsilon))\to \pi_i(X')$ is an isomorphism for $i<N-b$, and is an injection for $i=N-b$.
 Note that our methods can also be carried over to this 'topological' context; 
 they  yield (cf. part \ref{i3} of this remark) a comparison of the cohomology of $X'$ with the limit for $\eps\to 0$ of the cohomology of $s_X\ob(Z_\varepsilon)$  (if one wants to consider cohomology with integral coefficients here, then one can apply our methods for $\q$-coefficients and for $\z/l\z$-ones for all prime $l$ separately, and then combine the results obtained; see Theorem 5.2.16(ii) of \cite{dimca}). 
 This statement is somewhat weaker than the result cited; on the other hand the proof is much simpler than the one in ibid. (at least, it uses no stratified Morse theory). 
 
We can also 
 generalize our result in order to include the cases when 
 $U'$ 
 is  locally 'not quite a complete intersection'. Similarly to  \cite{gorma} (and to  Theorem 10.5(iv) of \cite{lyubez}), 
  this will decrease the highest degree in which we have a bijection (and an injection) for cohomology by a certain measure $c$ of the failure 
  for $U'$ to be  a LSTCI (i.e. by the number of 'extra equations' that are needed for defining $U'$ locally; see the Notation section).
In order to prove this statement only needs the corresponding generalization of Proposition \ref{ppt}(\ref{ilstci}) (cf. our argument above); the proof of loc. cit. can be easily extended to this setting. 

\item \label{i3}  We will discuss the difference between cohomological and homotopy formulations of Weak Lefschetz-type results in Remark \ref{rcfulaz} below.

Here we note the following:  if $\gi:Z\to X$ is an embedding of  paracompact topological spaces, $F$ is a 'topological' sheaf (of sets or abelian groups) on $X$, then 
\begin{equation}\label{eiim}
(\gi^*F)(Z)=\inli_{Z_0\supset Z,\ Z_0{\text{ is open in $X$}}}F(Z_0)
\end{equation}
(i.e., one does not need to sheafify here; see Corollary 1 in \S II3.3 of \cite{gode}).
Moreover, if $Z$ is compact and $X$ is a metric space then it suffices to take $Z_{\varepsilon}$ (the $\varepsilon$-neighbourhood of $Z$ in $X$) for all $\varepsilon>0$ instead of all possible $Z_0$ in (\ref{eiim}). Hence (cf. Proposition 5.6.6 of \cite{kanek}) Theorem \ref{tgara} yields indeed that a cohomological analogue of the result of \cite{gorma} (as discussed above) is fulfilled 'in the limit'. Besides,
considering the preimages of 'fat multiple hyperplane sections' (i.e., $s_X\ob(Z_\varepsilon)$) yields a method of computing $H^*(Z,\gi^*s_{X*}\znz_{X'})$ when $K=\com$. For a general $K$ one can consider certain henselizations instead; see Remark \ref{rhens} below.

\item \label{i4}  Taking $s_X=\id_X$ 
we immediately obtain the 'classical' Weak Lefschetz (for a LSTCI; for an arbitrary $X$ the corresponding bound on the degrees will decrease by the constant $c$ mentioned in part \ref{i2} of this remark). The only other way to prove this fact (for $p>0$) that is known to the author is to deduce it from Theorem 9.3 of \cite{lyubez}.


More generally, if $s_X$ is proper then proper base change  (see Proposition \ref{pbch}(I2)) yields that $\gi^*s_{X*}\znzr_{X'}\cong s_{Z*}\gi'^*\znzr_{X'}\cong  s_{Z*} \znzr_{Z'}$, i.e., we obtain an 'ordinary' Weak Lefschetz-type statement in this case also. The author was not able to find this result in literature. 

  Some  more interesting cases when Theorem \ref{tgara} yields that $H^*(X')\cong H^*(Z')$ are given by Corollary \ref{cpro} (where we consider $s_X$ that is proper only over an open neighbourhood of $Z$ in $X$) and by Theorem \ref{tbarth} below.

\item \label{i5}  
In contrast to \S II.1.2 of \cite{gorma}, we do not demand $s_X$ to be quasi-finite. 
Actually,  in \S II.1.1 of ibid. there is also a  version of a 'fat multiple hyperplane section weak Lefschetz' for a not necessarily quasi-finite 
$s_X$. Yet loc. cit. requires $X'$ to be smooth. This setting has certain advantages: 
instead of subtracting $e$ (i.e., the maximum of the dimensions of fibres of $s_X$) from $a-b$ to get the 
the bound on the homotopy degrees, in loc. cit. 
one only subtracts a certain measure $e'$ of the failure for $s_X$ to be semi-small.
Being more precise, if $\phi(k)$ for $k\ge 0$ denotes 
$\dim(\{{x\in X}:\ \dim s_X\ob(x)=k\})$
(here we set $\dim(\emptyset)=-\infty$), then  
one takes  
\begin{equation}\label{emin}
e'=\max_{k\ge 0}(2k+\phi(k)-a+ \min (\phi(k),b-1))+1-b \end{equation}
(this gives a better bound when $s_X$ is not equidimensional). 
Yet one has to pay a price for this refinement of the result; as can be seen from the example in \S II.8.4 of ibid., (both the cohomological and the homotopy version of) the statement with the improved bound  do not extend to the case of an arbitrary locally complete intersection $X'$. 

We note that Theorem \ref{tgara} implies the cohomological version of the statement mentioned (we have written it down  in Remark \ref{runiv}(2) below). 
Indeed,  for a smooth closed $Y\subset X'$ of (constant) codimension $c$
one can apply the Gysin distinguished triangle 
\begin{equation}\label{egysi5} v_*\znrz_{Y}(-c)[-2c]\to \znzr_{X}\to u_*\znrz_{X\setminus Y}
\end{equation} 
(see (\ref{eglu2})); here $v:Y\to X'$ and $u:X'\setminus Y\to X'$ are the corresponding embeddings). Hence instead of verifying the statement for $X'$ it suffices to verify it for $X'\setminus Y$ 'with the same $e'$' (i.e., up to the cohomological degree $a-b-e'$)  and to verify it for $Y$ up to the cohomological degree $a-b-e'-2c$ (certainly, there is nothing to check if $a-b-e'-2c<0$). Considering a smooth stratification of $X'$ 
such that the restriction of  $s_X$ to each stratum is equidimensional, one easily deduces from Theorem \ref{tgara} the result in question (by induction on the number of strata).

Also note that the Verdier dual to   (the $\znz$-module version of) Lemma III.7.4(1) of \cite{kw} 
 yields the following: $s_{U*}\znzr_{U'}\in \ds(U)^{p\ge a-e'}$ in the case when the maximum in (\ref{emin}) is attained for some value of $k$ such that $\phi(k)\ge b-1$.
Hence  one can prove the 'regular' version of Theorem \ref{tgara} in the case when this assumption on the maximum is fulfilled by slightly modifying the proof of loc. cit.  

Certainly, this 'regular' version of our main result immediately implies the corresponding versions of Corollary \ref{cpro} and of Theorem \ref{tbarth} below.

\item \label{i8} 
The only restriction on $U$ that we actually applied in the proof above is the restriction on its 'perverse' cohomological amplitude (note: the bound on the cohomological amplitude of $u_!$ from below is equivalent to the bound on the cohomological amplitude of $u_*$ from above by Verdier duality). In particular, one can consider a $U$ such that there exists an affine  bundle $p:U'\to U$ for which $U'$ is  affine (indeed, then we have that $p_*p^*\cong 1_{\ds(U)}$). 
Therefore the proof of Theorem 7.1.1 in \cite{lazar} yields the following: one can take for $U$ the complement of the zero locus 
of a section of an ample vector bundle (or rank $b$) over a projective $P$ (since we do not have to demand the transition maps of the corresponding bundle to be affine). 

Thus one can prove a certain extension of the Sommese's theorem; this result also implies the corresponding analogue of Corollary \ref{cpro} below. The latter statement seems to be quite new even in the case $K=\com$.



\item \label{i6}  The argument used in the proof of the Theorem can also be carried over to the relative context, i.e., for $X$ and everything else being schemes over a (more or less 'reasonable') 
 base scheme $S/\spe \zol$. Certainly, the conclusion will be formulated in terms of the 
 $t(S)$-cohomology of the corresponding total derived direct images. Then one can proceed to prove the corresponding (relative) analogues of Corollary \ref{cpro} and of
Theorem \ref{tbarth}(I) (below).

\item \label{i7}  
The author has also benefited from the proof of Theorem 6.1.1 of \cite{a1} (note that the result itself is not quite correct; cf. \cite{a07}).
Following loc. cit., we could have (directly) applied the Verdier duality instead of the  perverse $t$-structure 
in the proof of Theorem \ref{tgara}. 

Yet then we would (probably)  need the smoothness of $U'$. In fact, 
the perverse $t$-structure seems to be a very convenient tool for our purposes. 

\end{enumerate} 
\end{rema}

 Now we prove two easy statements, that often simplify the calculation of $H^*(Z,\gi^*s_{X*}\znzr_{X'})$.
 
 \begin{pr}\label{pinsub}

  Adopt the setting of Lemma \ref{lara}.

1. 
Suppose that there exists an open $O\subset X$ such that $Z\subset O$ and the restriction of $s_X$ to the preimage $O'$ of $O$ is proper. Then the base change transformation $B(s_X,\gi):\gi^*s_{X*} \implies s_{Z*}  \gi'^*$ is an isomorphism.


2. Let $\gi=t\circ v$, where $t:T\to X$ is a smooth morphism.
Then $\gi^*s_{X*}\znzr_{X'}\cong  v^*s_{T*}\znzr_{T'}$; here $s_T:T'\to T$ is the base change for $s_X$.

\end{pr}
\begin{proof}

1. Consider the commutative diagram
\begin{equation}\label{eqo} \begin{CD}
Z'@>{i_1'}>> O' @ >{i_2'}>> X'\\
@VV{s_Z}V
@VV{s_O}V @VV{s_X}V\\
Z@>{i_1}>> O @ >{i_2}>> X \end{CD}
\end{equation}
We have that $\gi^*=i^*_{1} i^*_{2}$. Proposition \ref{pbch}(II) yields that $B(s_X,\gi)$ is the composition of transformations $\gi^*s_{X*}\stackrel{B_1}{\implies} i^*_{1}   s_{O*} i'^*_{2}\stackrel{B_2}{\implies} s_{Z*}  i'^*_{1} i'^*_{2}=s_{Z*}  i'^*$. It remains to note that $B_1$ is an isomorphism by the
smooth base change theorem (applied to the right hand side square of (\ref{eqo})), whereas $B_2$ is an isomorphism by the  proper base change theorem (applied to the left hand side square of (\ref{eqo})); see   Proposition \ref{pbch}(I2).

2. Denote the base change of $t$ via $s_X$ by $t'$, and note that  $\gi^*=v^*t^*$. Then smooth base change yields:
$\gi^*s_{X*}\znzr_{X'}\cong v^* s_{T*} t'^*\znzr_{X'}=v^* s_{T*}\znzr_{T'}$.

\end{proof}

\begin{rema}\label{rthiso} 
Part 1 of the Proposition yields our (very simple) substitute of \S II5A of \cite{gorma} (that is sufficient for our applications below, and makes sense for any $p=\cha K$). 
 Note that loc. cit. relies on the Thom's first isotopy lemma whose proof is really long (see \S I1.5 of ibid.).
 
\end{rema}

Now we combine the previous result with Theorem \ref{tgara}. $\pi_1(-)$ will denote the \'etale fundamental group functor.

For simplicity in parts 2 and 3 of the Corollary below we will assume that $X'$ is connected. It follows (by part 1 of the Corollary) that  
$H^0(Z',\znz)\cong \znz$; hence $Z'$ is connected also, and we have no need to keep track of the corresponding base points.

\begin{coro}\label{cpro}
  Adopt the setting of Theorem \ref{tgara}  and Proposition \ref{pinsub}(1). Then the following statements are fulfilled.

1. $H^i(-,\znrz)(\gi')$ is bijective for all $i<a-e-b$, and is injective 
for $i=a-e-b$.

2. Suppose that $a-e-b \ge 1$. Then $\pi_1(\gi')$ is surjective. 

3. Assume that $a-e-b \ge 2$. Then the kernel of $\pi_1(\gi')$ has no factors isomorphic to $\zlz$ (for all prime $l\neq p$).


\end{coro}
\begin{proof}
1. First note that $H^i(X',\znzr_{X'})\cong H^i(\spe K,x'_*\znrz_{X'})=H^i(\spe K,x_*s_{X*}\znrz_{X'})$ and $H^i(Z',\znzr_{Z'})\cong H^i(\spe K,z'_*\znzr_{Z'})= H^i(\spe K, x_*s_{X*}\gi'_*\gi'^* \znzr_{X'})$ (by Proposition \ref{pds}(\ref{icohh})). By Proposition \ref{pbch}(I1) we obtain that $H^i(-,\znrz)(\gi')$ is given by applying $H^i(\spe K,x_*(-))$ to the composition of $\gi_*(B(s_X,\gi)(\znrz_{X'}))$ with $\gff_r$.

Hence the statement can be obtained immediately by combining  Theorem \ref{tgara}  with the first part of 
Proposition \ref{pinsub}.

2. 
Consider a finite connected etale cover $X''\to X'$, 
and the corresponding cartesian square 
\begin{equation}\label{eqfund} \begin{CD}
Z''@>{\gi''}>>  X''\\
@VV{}V
@VV{}V \\
Z'@ >{\gi'}>> X' \end{CD}
\end{equation}
Now, we can apply assertion 1 for $\gi''$ instead of $\gi'$; hence $H^0(-,\znz)(\gi'')$ is bijective, i.e., $Z''$ is connected also. 

Next we recall 
that finite connected etale covers of any variety $Y$ are given by open subgroups of $\pi_1(Y)$,
and that the degree of such a cover equals the index of the corresponding subgroup (in $\pi_1(Y)$).
Hence for the subgroup $H$ of $\pi_1(X')$ corresponding to $X''$, we obtain a bijection of cosets $\pi_1(Z')/(\pi_1(\gi')\ob (H))\to \pi_1(X')/H$. 

Considering $H$ running through
a projective system of open subgroups 
 of $\pi_1(X')$ such that their limit is  $\ns$, we obtain 
the result.

3. We use the notation and 
the 
statements mentioned in the proof of the previous assertion. 
Moreover, applying assertion 1 to $X''$ instead of $X'$ we obtain that $H^1(-,\znz)(\gi'')$ is bijective also. 

Now  recall that $H^1(Y,\znz)\cong \homm (\pi_1(Y),\znz)$. Hence taking a 'small enough' $H$ (i.e., letting it to run through a system of open subgroups  of $\pi_1(X')$ such that their limit is  $\ns$) again ,we obtain the result. 

\end{proof}

\begin{rema}\label{rpi1}
We conjecture that in the setting of Corollary \ref{cpro}(3) the homomorphism $\pi_1(\gi')$ is actually bijective. 
Possibly, one can prove this conjecture in general by combining our methods with those of \cite{cutpi1} (somehow).

\end{rema}

\section{Applications: certain Barth-type theorems} \label{sappl}

It turns out that 
Corollary \ref{cpro}(1) allows to study the cohomology of the intersection of a (closed LSTCI) subvariety of $(\p^N)^q$ with the diagonal (as well as the cohomology of  the preimages of the diagonal with respect to proper morphisms whose target is $(\p^N)^q$). As a consequence, we can calculate the lower cohomology of (the preimages of) subvarieties of $\p^N$ (of small codimension); so we extend the seminal results of Barth and others 
 to the 
 case of arbitrary characteristic ($\neq l$).



Our exposition closely follows the one of \S9 of \cite{fulaz} (see also \S3.5.B of \cite{lazar}). The main distinctions are due to the fact that we consider the cohomology groups of varieties instead of the homotopy ones. In particular, this makes the proof of our Theorem \ref{tbarth}(II1,2) quite different and somewhat more complicated than the proofs of Corollaries 9.7 and 9.8 of  \cite{fulaz}. On the other hand, in our Theorem \ref{tbarth}(II2) we do not require $v$ to be finite (in contrast to Corollary 9.8 of \cite{fulaz}).


We will need some notation. We fix some $l,n,r$, and also some $N>0$. 

We recall a construction of  certain locally trivial $G_m^{s}$-bundles (for  $s\ge 0$); 
it is described in more detail in \S3 of ibid. (for the case $q=2$).

For any $q>0$
we note that the natural projection $a_q:(\af^{N+1}\setminus \ns)^q\to (\p^N)^q$ factors through the quotient $V_q$ of $(\af^{N+1}\setminus \ns)^q$ by the diagonal action of $G_m$ on $\af^{Nq+q}\setminus\ns\supset (\af^{N+1}\setminus \ns)^q$. So, we obtain a (locally trivial) $G_m^{q-1}$-bundle $p_q:V_q\to (\p^N)^q$, and a 
 $G_m$-bundle $b_q:(\af^{N+1}\setminus \ns)^q\to V_q$, whereas $V_q$ is an open subvariety of $\p^{Nq+q-1}$.

For a morphism of varieties $g:Y\to (\p^N)^q$ 
we will denote by $g_q:Y_q\to V_q$ (resp. by $g'_q:Y'_q\to (\af^{N+1}\setminus \ns)^q$) the base change of $g$  along $p_q$ (resp.  along $a_q$). 


\begin{theo}\label{tbarth}

I Let $q>1$.

1. 
 Let 
 $g:Y\to (\p^N)^q$ be a proper morphism (of varieties) 
 of relative dimension $\le e$, where $Y$ is a LSTCI (see the Notation section) of dimension $a$. 
  Denote the diagonal of $(\p^N)^d$ by $\Delta(\cong \p^N)$. Then 
 there is a 
 natural morphism $c: g\ob(\Delta)\to Y_q$ 
such that $H^i(-,\znrz)(c)$ is bijective for $i<a-e-qN+N$ and is injective for $i=a-e-qN+N$.

2. 
In the setting of the previous assertion 
assume that $Y=(\prod Y_j)_{red}$ (see the Notation section) for some proper morphisms $g_j:Y_j\to \p^N$ ($1\le j\le q$), and  $g$ is the restriction of $\prod g_j$ to $Y$ (if $K$ is perfect, then we always have $g=\prod g_j$ and $Y=\prod Y_j$).
Then $Y_q'\cong \prod Y'_{j1,red}$, and for the base change $c': (g\ob(\Delta))_1'\to \prod Y'_{j1,red}$
of the morphism $c$ (that is given by the previous assertion) to $Y_q'$ 
we have that $H^i(-,\znrz)(c')$ is bijective for $i<a-e-qN+N$ and is injective for $i=a-e-qN+N$.




II Let $t:T\to \p^N$ be  a closed embedding, where $T$ is a LSTCI of dimension $d$. Then 
the following statements are valid.

1. $H^i(-,\znrz)(t)$ is bijective for $i\le 2d-N$ and is injective for $i=2d-N+1$;
the same is true for the corresponding  $G_m$-bundle morphism $t_1':T_1'\to \af^{N+1}\setminus \ns$.

2. Let $v:V\to\p^N$ be a proper morphism of relative dimension $\le b$, where $V$ is a LSTCI of dimension $u$. Then for the morphism $w: v\ob (T)\to V$ we have that $H^i(-,\znrz)(w)$ is bijective for $i\le \min (d+u-b-N-1, 2d-N)$ and is injective for $i=\min (d+u-b-N,2d-N+1)$. 
Besides, the same is true for 
$w_1'$ being the base change of $w$ to $V_1'$.


\end{theo}
\begin{proof}
1. The diagonal embedding $\af^{N+1}\setminus \ns \to (\af^{N+1}\setminus \ns)^q$ yields (after the factorization by the 
diagonal action of $G_m$) a subvariety $L_q\subset V_q$ such that the restriction $p_L$ of $p_q$ onto $L_q$ gives an isomorphism $L_q\to \Delta$. Note also that $L_q$ is a closed subvariety in  $\p^{Nq+q-1}\supset V_q$.

We denote the embedding $L_q\to \p^{Nq+q-1}$ (resp. $
\Delta\to V_q$)  by $\gi$ (resp. by $\gi_1$). Take for $c$ the base change of $\gi_1$ along  $g_q$, and denote the composite morphism $Y_q\to V_q\to \p^{Nq+q-1}$ by $s$.

Now we apply Corollary \ref{cpro}(1) as follows: we take $s_X=s$, $Z=L_q$, $O=V_q$. Since $g_q\ob(L_q)\cong g\ob(\Delta)$, we obtain the result.

2. The first part of the assertion is obvious. 
Next we note that $\Delta\cong \p^N$, 
 and the base change of $c$  to $Y_q'$ does yield $c':(\prod_{\p^N} Y_j)_{1,red}'\to (\prod Y'_{j1})_{red}$. 

The cohomological part of our statement follows immediately from the previous assertion (applied for all $r\in\z$) 
together with  Lemma \ref{lgmbundle} 
 (below). 

II By Proposition \ref{ppt}(\ref{ihsss}), we can assume that $K$ is algebraically closed. In this case we  can (and will) set $r=0$ (see Proposition \ref{pds}(\ref{ikun})). We denote the functor $H^*(-,\znz)$ by $H^*$.

1. 
Lemma \ref{lgmbundle} yields that it suffices to prove 
the second part of the assertion. Since the \'etale cohomology of $\af^{N+1}\setminus\ns(K)$ vanishes in all degrees between $1$ and $2N$,
to this end it suffices to verify that $H^i(t_1')$ is an isomorphism for all $i\le 2d-N$.

Now in the notation of assertion I.2 we take $q=2$ and 
 $g_1=g_2=t$. For the corresponding $Y=T\times T$ (that also equals $(T\times T)_{red})$ we have $a=2d$, $e=0$.

Hence assertion I.2 yields that for the corresponding $c':T'_1\to T'_1\times T'_1$ we have that $H^i(c')$ is bijective for $i<2d-N$, and is injective for $i=2d-N$.

Now we note that the projection $pr:T_1'\times T_1'\to T_1'$ (via the first factor) splits $c'$; hence $H^*(pr)$ yields  an isomorphism up to degree $2d-N-1$ and a surjection in  degree $2d-N$. 

Now we 
rewrite the result obtained in terms of $D(Sh^{et}(\spe K, \znz-\modd))\cong D(\znz-\modd)$ and of the total derived functors $RH(-)\cong p_*$ (where $p:T_1'\to \pt=\spe K$ is the structure morphism; see Proposition \ref{pds}(\ref{ik},\ref{icoh})). 
 Choose some morphism $s:\pt \to T_1'$. Since $p\circ s=\id_{\pt}$, we obtain a splitting
 $RH(T_1')\cong \znz(=RH(\pt))\bigoplus D$, where $D$ is the  'image' of the idempotent endomorphism $(\id_{RH(T_1')}- RH(s\circ p))\in D(\znz-\modd)(RH(T_1'), RH(T_1'))$. Now, the Kunneth formula (see Proposition \ref{pds}(\ref{ikun})) yields that $RH(pr)$ can be described as the tensor product of the split morphism $\znz\to \znz \bigoplus D$ by $RH(T_1')\cong \znz\bigoplus D$. Hence $H^i(pr)$ is injective for any $i\ge 0$. 
 Moreover, we obtain that the cohomology of $D$ is concentrated in degrees $>2d-N$, since $D$ is a direct summand of $\co(RH(pr))\cong D\bigoplus D\otimes D$; this concludes the proof. 

2. Again, it suffices to prove the second part of the assertion. We set $q=2$ (again) and take $g_1=t$, $g_2=v$.
For the corresponding $Y=T\times V$ we have $a=d+u$, $e=b$.

Hence assertion I2 yields that for the corresponding $c':(v\ob (T))_{1}'\to V_1'\times T_1'$ we have that $H^i(c')$ is bijective for $i<d+u-b-N$ and is injective for $i=d+u-b-N$.

Now, assertion II1 (combined with Proposition \ref{pds}(\ref{ikun})) 
yields that for the projection $pr:V_1'\times T_1'\to V_1'$ we have the following: $H^i(pr)$ is a bijection for $i\le 2d-N$ and is an injection for $i=2d-N+1$. Since $H^i(w_1')=H^i(c')\circ H^i(pr)$, we obtain the result.  

\end{proof}

\begin{rema}\label{rlyu}

1. 
The corresponding Leray spectral sequence yields a relation of the cohomology of $Y_q$ with the one of $Y$; for $q=2$ one obtains  a Gysin long exact sequence (see  (\ref{efrgys})).

2. Applying our cohomological results to $H^*(-,\znz(1))$, one immediately obtains certain statements about the Picard and Brauer groups of varieties mentioned (if the corresponding bound on the cohomology degrees where we have an isomorphism 
is not too small). One can also apply our results to get some information on the number of points of varieties over finite fields. 

The author does not know (at the moment) whether the corresponding statements are interesting (and whether these results for Brauer and Picard groups are new; cf. \S11 of \cite{lyubez} for some results on the Picard groups). 


3. 
Note 
 that in contrast to \S9 of \cite{fulaz} we do not demand $f$ to be finite; cf. Remark \ref{rara}(\ref{i5}) above.

4. Certainly, our results can be extended to the case when  the variety $Y$ considered above (this includes $Y=T\times T$ and $Y=T\times V$ in the proofs of assertions II1 and II2 of the theorem) is not ('quite') a LSTCI; the corresponding degree bounds should be decreased by a certain measure $c$ of their failure to be so (cf. Remark \ref{rara}(\ref{i2})). 
Also, it is not really necessary to assume that $Y$ is equidimensional.
A certain extension of our Theorem \ref{tbarth}(II1) of this sort along with several other results (on the cohomology of subvarieties of $\p^N$ of low codimension) 
can be found in Theorem 10.5 of \cite{lyubez}.

5. In particular, one can easily prove the following statement: if $q>1$ and $T_j\subset \p^N$ are closed  LSTCI of dimensions $d_j$ for $1\le j\le q$, then  $H^i(\cap_{1\le j\le q} T_j,\znrz)\cong H^i(\p^N,\znrz)$ for $i\le D=\sum_{j=1}^q d_j-(q-1)N$. 
 Indeed, we can assume that $d_1$ is the largest of $d_j$; take $T=T_1$ and $V=\cap_{2\le j\le q}T_j$. Then $H^i(T,\znrz)\cong H^i(\p^N,\znrz)$ for $j\le 2d_1-N$ 
 by part II1 of our theorem. 
 In the case when $T_j,\ 2\le j\le q$, intersect properly (i.e. the dimension of their intersection is $\sum_{j=2}^q d_j-(q-2)N=D+N-d_1$) it remains to apply part II2 of our theorem; otherwise one can apply the 'defect' version of this result mentioned above.

Thus we are able to extend the statement mentioned after Remark 9.9 of \cite{fulaz} to the case of an arbitrary $K$; it seems that this assertion does not follow from the results of \cite{lyubez}.



6. Part 
II of our theorem can easily be 'relativized'. 
Let $S$ be any 
 $\spe \zol$-scheme. We will consider schemes $w:W\to S$ such that the functor $w_*$ has finite cohomological dimension.
  This means the following: we consider the \'etale sheaves $H^i_S(W, K)\in \obj Sh^{et}_{\znz-\modd}(S)$ being the 
cohomology of $w_*K$, $K\in \obj Sh^{et}_{\znz-\modd}$ (i.e. we consider the canonical $t$-structure on $\ds(S)$ and not the perverse one). 
 We will
 call such a 
scheme $W/S$ a (relative)
 LSTCI if for each Zariski (or \'etale) point $f:\spe K\to S$ the fibre $W_K$ of $W$ over $\spe K$ or the associated reduced scheme $W_{K, red}$ is a LCSTI variety over $K$ of dimension that does not depend on the choice of 
 $f$. We will use similar conventions for dimensions and relative dimensions (i.e. all of them are 'measured at points of $S$'). Note here: one can easily define natural 'global' analogues of these notions (at least) in the case when $S$ is excellent noetherian of finite Krull dimension.
 
  Let $t:T\to \p^N$ be  a closed embedding, where $T$ is a (relative) LSTCI of dimension $d$ (in our 'relative' sense). Then 
the following statements are valid.

(i)  $H_S^i(-,\znrz)(t)$ is bijective for $i\le 2d-N$ and is injective for $i=2d-N+1$. 

(ii) Let $v:V\to\p^N$ be a proper morphism of relative dimension $\le b$, where $V/S$ is a LSTCI of dimension $u$. Then for the morphism $w: v\ob (T)\to V$ we have that $H^i(-,\znrz)(w)$ is bijective for $i\le \min (d+u-b-N-1, 2d-N)$ and is injective for $i=\min (d+u-b-N,2d-N+1)$.

Indeed, recall that in order to compare the cohomology sheaves in question it suffices to calculate their stalks 
at (\'etale or Zariski) points of $S$. Proper base change (see Proposition \ref{pbch}(I2); it can be applied in this more general setting by Theorem 4.3.1 of \cite{sga417}) yields that the stalks in question are the \'etale cohomology 
groups of the corresponding varieties 
 over $\spe K$. Hence parts II1 and II2 of Theorem \ref{tbarth} imply their relative analogues indeed.

7. 
Possibly there exist some other (locally trivial)  $G$-bundle constructions (where $G$ is some algebraic group) similar to the  $G_m$-bundles considered above, such that our methods can be used for them in order to obtain certain Barth-type 
results. 

\end{rema}


In order to conclude the proof 
of Theorem \ref{tbarth}, it remains to prove the following statement.

\begin{lem}\label{lgmbundle} Let  $f':A'\to B'$ be a morphism of locally trivial $G_m$-bundles
over the base $f:A\to B$; 
fix a $j\in\z$.
 Then we have that $H^i(-,\znrz)(f)$ is bijective for all $i<j$ and is injective for $i=j$ and all $r\in \z$ if and only if the same is true for $H^i(-,\znrz)(f')$.

\end{lem}
\begin{proof}
We prove the statement by induction in $j$.

 If $j<0$, then the statement is vacuous. Now assume that it is fulfilled for $j=s-1$, and that
 $H^i(-,\znrpz) (f)$ and $H^i(-,\znrpz)(f')$ are bijective for all $i<s-1$ and are injective for $i=s-1$ (for some $s\ge 0$, and for all $r'\in \z$). We should verify (for some fixed $r\in \z$) that $H^{s-1}(-,\znrz)(f)$ is surjective and $H^{s}(-,\znrz)(f)$ is injective if and only if the same is true for $f'$.

Now assume that for any 
$G_m$-bundle $b:X'\to X$  there exists a 
certain (Gysin) long exact sequence  
\begin{equation}\label{efrgys} \begin{aligned} 
 H^{s-1}(X',\znrz)\to H^{s-2}(X,\znroz)\to H^s(X,\znrz) \to \\
 H^s(X',\znrz)\to H^{s-1}(X,\znroz) \to  H^{s+1}(X,\znrz) 
  \end{aligned}
\end{equation} 
and that this sequence is functorial with respect to $b$.

Then the five lemma yields the following: if $H^{s-1}(-,\znrz)(f)$ is bijective and $H^{s}(-,\znrz)(f)$ is injective, then $H^{s-1}(-,\znrz)(f')$ is bijective; if  $H^{s-1}(-,\znrz)(f')$ is bijective, then $H^{s-1}(-,\znrz)(f)$ is bijective also. Next, the four lemma on monomorphisms yields the following: if $H^{s}(-,\znrz) (f)$ is injective  then $H^{s}(-,\znrz)(f')$ is injective also;  if $H^{s-1}(-,\znrz)(f')$ is bijective and $H^{s}(-,\znrz)(f')$ is injective, then $H^{s}(-,\znrz)(f)$ is injective.

Hence we obtain the inductive assumption for $j=s$.

Now it remains to verify the existence and the functoriality of (\ref{efrgys}). 
We note that $X'$ can be presented as the complement to $X$ for the line bundle $X''=X\times \af^1/G_m$ over $X$ (we consider the diagonal action of $G_m$ on the product and the zero section embedding $X\to X''$). Hence
it remains to apply Corollary 1.5 of \cite{sga5e7}. 


\end{proof}

\begin{rema}\label{rcfulaz}
Now we say more on the relation of our results with those of \cite{fulaz} (where for 
$K=\com$ the homotopy analogues of all of the results of this section were proved), 
and also explain why our results are formulated in the way they are, and how could they be modified.

1. \S9 of \cite{fulaz} relies on the results of \cite{gorma}. So loc. cit. has three features that distinguish it from our 
Theorem \ref{tbarth}:

(i) it treats only complex varieties;

(ii) all the corresponding morphisms to $\p^N$ and $(\p^N)^q$ are required to be finite;

(iii) the results are formulated in terms of homotopy instead of cohomology.

Certainly, (i) is a serious disadvantage of \cite{fulaz}; the restriction (ii) probably could be removed via generalizing Theorem II1.2 of \cite{gorma}.

2. Now we discuss the difference between treating cohomology and homotopy.

Certainly, cohomology is much easier to compute. On the other hand, homotopy groups have certain theoretical advantages. 
One of them is that the relation between the (higher) homotopy groups of a (locally trivial) $G_m(\com)$-bundle
$X'\to X$ with those of $X$ are much easier to use than  (\ref{efrgys}). Using this, it was proved (in Theorem 9.2 of loc. cit., whose setting corresponds to that of our Theorem \ref{tbarth}(I1) for $q=2$) that the homotopy groups of $Y$ are isomorphic to those of $g\ob(\Delta)$ in degrees between $2$ and $a-N$ (note that $e=0$ in loc. cit.). Moreover, in the setting of our assertion I2 (with $q=2$) 
this result reformulates as an isomorphism (in the corresponding degrees) of relative homotopy groups for the pairs $(Y_1,Y_1\times_{\p^N} Y_2)$ and $(\p^N,Y_2)$.  The homotopical analogues of our assertions II(2,1) follow easily. 

So, it would be interesting to obtain certain homotopy analogues of the results above. Now we discuss  possible ways to do this. First we note that even for complex varieties it is impossible to recover the 'topological' homotopy groups staying inside the category of algebraic varieties. Therefore, one has to consider the homotopy groups of the corresponding \'etale homotopy types.

To the knowledge of the author, the existing homotopy-theoretic technique does not allow to modify the proof of Theorem \ref{tgara} so that it would yield information on (\'etale) homotopy groups directly. For this reason, in order to obtain the homotopy analogues desired one would probably have to apply the results above, and then extract the information on homotopy from the one on cohomology using 
 \'etale analogues of the Hurewicz theorem. 

An obvious obstruction to do so is the (possible) non-triviality of $\pi_1$ of the varieties in question. The author hopes that   Corollary \ref{cpro}(2,3) or the methods of its proof can help here (cf. also Remark \ref{rpi1}).
 

3. Now we describe a possible 'cohomological' way to get rid from the $G_m$-bundles in (the formulation of) Theorem \ref{tbarth}(I). 
 As can be shown by simple examples, there is no 'easy' way to do this. Yet 
 the author suspects that (similarly to the case of a pullback square of topological spaces, as studied by Eilenberg and Moore)
  in the setting of Theorem \ref{tbarth}(I2)
 the differential graded algebra that computes the cohomology of $\prod_{\p^N} Y_i$ is 
  quasi-isomorphic up to degree $a-e-qN+N$ to the tensor product of the corresponding algebras for $Y_i$ over the one for $\p^N$ (maybe, only when $K$ is separably closed). 
 Note that for the corresponding $G_m$-bundles (i.e., for $(\prod_{\p^N} Y_i)_1'$ and $\prod Y'_{i1}$) this fact is given by  loc. cit.; then one probably  has to 'descend' using an (induction) argument similar to the one in the proof of Lemma \ref{lgmbundle}.

  


\end{rema}

\section{Some other remarks: henselizations as 'small neighbourhoods'}\label{srem}

\begin{rema}\label{runiv}
1. As we have already noted, Theorem \ref{tgara} generalizes the 'classical' Weak Lefschetz. Moreover, the author knows no way to deduce loc. cit. from the 'ordinary' Weak Lefschetz (or even to reduce it to the case when $s_X$ is proper). 

Yet we note that the 'regular version' of loc. cit. mentioned in Remark \ref{rara}(\ref{i5})  can be reduced to its rather special modification. 

2. First we formulate the 'regular version' here. 

In the setting of Theorem \ref{tgara} we assume that $X'$ is regular. 
For any $k\ge 0$ denote by
$\phi(k)$ the  dimension of 
$\{{x\in X}:\ \dim s_X\ob(x)=k\}$
(here we set $\dim(\emptyset)=-\infty$).  Then  
  for $\beta =2a-1 -\max_{k\ge 0}(2k+\phi(k)-a+ \min (\phi(k),b-1))$ the homomorphism $H^i(X,-)(\gff_r))$  is a bijection if $i< \beta$, and is an injection if $i=\beta$.

3. Now we verify briefly that the general case of the statement above follows from its '\'etale-local version' corresponding to $X'$ being a geometric point of $X$ (considered as an ind-\'etale scheme over  the corresponding 
closed subvariety of $X$).

Indeed, using the Gysin distinguished triangle (\ref{egysi5})  
 and passing to the limit we obtain the following: the general case of the statement follows from its analogue for $X'$ being 
a (regular) variety over a Zariski point $P\in X$. Moreover, applying \'etale descent one obtains:
it suffices to prove the statement for $X'$ being 
a variety over a geometric point $P_{geom}/P$.

 As in the proof of Theorem \ref{tgara}, we assume that $n=1$. Denoting the corresponding morphism $X'\to P_{geom}$ by $s_{geom}$, and denoting the  morphism $P_{geom}\to X$ by
$j_{geom}$, we obtain the following: we should compare  $A=H^*(X, j_{geom*}s_{geom*}\znzr_{X'})$ with $B=H^*(Z,\gi^*j_{geom*}s_{geom*}\znzr_{X'})$ (cf. the proof of Theorem \ref{tgara}). Now, the (small) \'etale site of  $P_{geom}$ is trivial, and $s_{geom*}\znzr_{X'}\in \ds(P_{geom})^{p\ge 0}$; hence $s_{geom*}\znzr_{X'}$ splits as a direct sum of $\znz_{P_{geom}}[r_i]$ for some $r_i\le 0$. Hence instead of $s_{geom*}\znzr_{X'}$ both in A and B one can put $\znz_{P_{geom}}$ (that is isomorphic to $\znzr_{P_{geom}}$ for any $r\in \z$); this finishes the reduction in question. 

4. Besides, an obvious \'etale descent reasoning allows to reduce Theorem \ref{tgara} (for $s_X:X'\to X$, $X'$ can be singular) to its analogue for the strict henselizations of Zariski points of $X'$.

Note that that the reductions described could be called 'motivic'.

\end{rema}

 This observation fits nicely with the 'henselization' methods for calculating $\gi^*$, that we  
will consider in the next remark.

\begin{rema}\label{rhens}
Now we describe some 
methods for studying the cohomology of $\gi^*s_{X*}\znrz_{X'}$ using certain henselizations.
They hardly have any computational value; their advantage is that they are somewhat 'motivic' and can be applied when no properness assumptions on $s_X$ (cf. Corollary \ref{cpro}) are fulfilled. 

By the main result of \cite{gabaffbch} and \cite{rhub}, if $h:A\to B'$ is a henselian pair (of affine schemes) and
$C'$ is a complex of \'etale sheaves over $B'$, then $H^*(B', C')\cong H^*(A, h^*(C'))$. Therefore if we decompose a closed embedding $i:A\to B$ of affine schemes as $A\stackrel{h}{\to} B'\stackrel{e}{\to} B$, where $B$ is the henselization of $A$ in $B$ and $e$ is the corresponding pro-\'etale morphism, then for a complex $C$  of \'etale sheaves over $B$ we have $H^*(A,i^*C)\cong 
H^*(B', e^*C)$. Now, since $e$ is a pro-\'etale morphism, $e^*$ is just the corresponding 'restriction' functor, and it 'commutes with base change'.

Thus, henselizations in $X$ could be thought about as of algebraic analogues of  'very small' neighbourhoods of  (closed) submanifolds $Z$ of a manifold $X$ (the latter can be used in order to compute the 'topological' functor $\gi^*$; see Remark \ref{rara}(\ref{i3})). Unfortunately, 'nice' henselizations exist only for affine schemes, whereas
$X$ and $Z$ are proper. We describe some possible methods for overcoming this difficulty.


1.  
 The first way 
 is to  choose an \'etale affine hypercovering $X.\to X$ (certainly, $X.$ can also be a Zariski or a Nisnevich hypercovering; it could 
 be a Cech hypercovering). Then for the corresponding morphisms $c_j:X_j\to X$ and their restrictions $c_j':c_j\ob(Z)\to Z$ (for $j\ge 0$) we would have a spectral sequence $H^i(c_j\ob(Z),{c'}_j^*\gi^*s_{X*}\znrz_{X'})\implies H^{i+j}(Z,\gi^*s_{X*}\znrz_{X'})$, whereas smooth base change yields that 
the $E_1$-terms of this spectral sequence are the corresponding cohomology of the henselizations of $c_j\ob(Z)$ in $X_i$ (for  non-connected schemes the henselizations are defined componentwisely).

2. Another way to reduce the computation in question to a one for affine schemes is to apply Jouanolou's device.

In the setting of Theorem \ref{tgara} let $b_X:\hat{X}\to X$ be an 
 affine vector bundle torsor such that $\hat{X}$ is affine; 
we denote the pullback of $b_X$ to $X'$ and $Z$ by $b_{X'}:\hat{X'}\to X'$ and $b_Z:\hat{Z}\to Z$, respectively; we will also use similar notation for morphisms.
Since the cohomology of the fibres of $b_Z$ over geometric points is trivial, we have that $b_{Z*}b_Z^*\cong 1_{\ds(Z)}$.
For $C\in \obj\ds(X)$ we obtain the following: $H^*(Z,\gi^*C)\cong H^*(\hat{Z},b_Z^*\gi^*C)\cong H^*(\hat{Z}, \hat{\gi}^*b_X^*C)$ (we used smooth base change in the last isomorphism). Now, $\hat{\gi}$ is a closed embedding of affine varieties, and one obtains that the cohomology groups in question 
are isomorphic to $H^*(\hat{Z}_h,\hat{z}_h^*b_X^*C)$, 
where $\hat{Z}_h$ is the henselization of $\hat{Z}$ in $\hat{X}$, $\hat{z}_h:\hat{Z}_h\to \hat{X}$ is the corresponding pro-\'etale morphisms. 

Then we obtain the following: $$  \begin{aligned} H^*(Z,\gi^*s_{X*}\znzr_{X'}) \cong
H^*(\hat{Z}_h,\hat{z}_h^*b_X^*s_{X*}\znzr_{X'}) \\
\cong H^*(\hat{Z}_h,s_{\hat{Z}_h}  \hat{z'}_h^*{b'}_X^*\znzr_{X'})
\cong H^*(\hat{Z}'_h,\znrz),\end{aligned}$$
 where $\hat{Z}'_h\stackrel{\hat{z'}_h}{\to}\hat{X}\stackrel{b'_X}{\to} X$ is obtained by base change via ${s_X}$ from 
$\hat{Z}_h\stackrel{\hat{z}_h}{\to}\hat{X}\stackrel{b_h}{\to} X$ (we apply smooth base change again).

3. 
One more  application of the smooth base change (in the setting of Theorem \ref{tgara})
can be obtained by applying Proposition \ref{pinsub}(2).
In particular, if $Z$ is a (multiple) hyperplane section of $X$ (corresponding to some embedding of $X$ into $\p^N$ for some $N>0$) then 
it can be presented as a member of a 'smooth family' of (multiple) hyperplane sections of $X$. This means: 
for some variety $B$ and some closed $T\subset B\times X$
 the projection $t:T\to X$ is a smooth morphism 
and  it restricts to an isomorphism $b\op(p)\to Z$,
where  $p$ is some closed point  $B$, $b$ is the projection $T\to B$
(whereas all closed fibres of $b$ yield multiple hyperplane sections of $X$ via $t$). 
Now, one can assume that $B$ is affine; denote by $p_h$ the henselization of $p$ in $B$, and denote by $Z_h$ the base change of $p_h$ with respect to $b$. 
Then (a very simple case of) Corollary 1 
of \cite{gabaffbch} yields the following (in the notation of  Proposition \ref{pinsub}(2)): $H^*(Z,\gi^*s_{X*}\znzr_{X'})\cong  H^*(Z,v^*s_{T*}\znzr_{T'})\cong H^*(Z_h, v_h^*s_{T*}\znzr_{T'})\cong H^*(T'_h, \znrz)$, where $v^*$ is the isomorphism $Z\to b\op(p)$, $v_h$ is the corresponding pro-\'etale morphism $Z_h\to T$, and  $T'_h$ is the base change of $Z_h/T$ with respect to 
$s_T$. 

 For $K=\com$ it can be easily verified that one can take an 'infinitely small ball' around $p$ instead of $p_h$ (cf. Remark \ref{rara}(\ref{i3})). So, our observation above  could be thought of as being one more version of an (algebraic) 'fat multiple hyperplane section'
approach to Weak Lefschetz-type questions.  
 It is also closely related with the  classical yoga of general hyperplane sections (of not necessarily projective varieties; cf. Lemma 3.3 of \cite{bei87}), especially in the case when $X$ is smooth and $s_X$ is an open embedding. 

For all of the methods listed in this remark it seems interesting to consider $X'$ that  runs through  geometric points of $X$ (cf. Remark \ref{runiv}).

4. The current paper grew our from an attempt to prove a certain weak Lefschetz-type result for $\znz$-motivic cohomology (at least, of complex varieties).
The 
idea was to join Theorem II1.2 of \cite{gorma} together with a formula that relates the motivic cohomology with the \'etale one (see Corollary 7.5.2(2) 
 of \cite{bws}; this is a more or less simple consequence of the Beilinson-Lichtenbaum conjecture) in order to study torsion motivic cohomology. Unfortunately, it turned out that in order to realize this program (which certainly becomes more realistic if one replaces Theorem II1.2 of \cite{gorma} with our Theorem \ref{tgara}) one needs a certain 'henselian model' for $Z$ in $X$ (a sort of 'motivic tubular neighbourhood'; cf. \cite{lemotubnei}) that is ind-\'etale 
 over $X$. Parts 2 and 3 of this remark do not yield such a model, whereas in part 1 we obtain a 'complex' of ind-\'etale $X$-schemes. 
One could also try to consider some $G_m$-bundles over the varieties in question (similarly to \S\ref{sappl}; cf. also \S10 of \cite{lyubez}). 
 
 Yet the author has some ideas for overcoming this difficulty (and to prove at least that the lower motivic cohomology of $X$ is isomorphic to the one of the scheme $\hat{Z}_h$ that was considered in part 2 of this remark).  It also seems to make sense to pass to the limit with respect to Nisnevich hypercoverings $X.\to X$ (in part 1 of this remark), and use the fact that $\gi^*_{Nis}(\znz(r)_{Nis,X})\cong \znz(r)_{Nis,Z}$ (here $\znz(r)_{Nis,-}$ is a complex of Nisnevich sheaves that computes the corresponding motivic cohomology; this isomorphism is given by a certain 'rigidity' argument). 
 

\end{rema}

\end{document}